\documentclass[12pt]{iopart}

\font\Bbb=msbm10 scaled 1000
\newcommand{\N}{\mbox{\Bbb N}}
\newcommand{\R}{\mbox{\Bbb R}}
\newcommand{\CC}{\mbox{\Bbb C}}

\newcommand{\fm}{{\phantom{-}}}

\begin{document}

\title[Spectral properties of the period-doubling operator]{Spectral properties of the period-doubling operator}

\author{V P Varin}

\address{Keldysh Institute of Applied Mathematics, Miusskaja sq. 4, 125047, Moscow, Russia}

\ead{varin@keldysh.ru}

\begin{abstract}

We compute the spectrum of the Feigenbaum period-doubling operator
in the space of bounded analytical functions in an ellipse.
The spectral properties of the period-doubling operator
in this space are not the same as in the space of even analytical functions.
In particular, it was found that the dimension of the unstable
manifold is not one (Feigenbaum's conjecture), but three.
We analyze several articles devoted to this problem
and compare different approaches and algorithms.

\end{abstract}
\maketitle

\section{Introduction}

Since the publication of the seminal works by Feigenbaum \cite{Fg1,Fg2},
hundreds of studies were devoted to this very interesting and still
expanding field of research. The author of this study was familiar
for many years with the subject, but in the most general terms.
One of the subjects that the author specializes in is
the development of efficient numerical and
symbolic algorithms for solving various mathematical problems.
In the course of testing one of such algorithms for solution of
functional equations, the Feigenbaum universality equation appeared to be a very
convenient model problem. The fundamental constants associated
with this equation are computed to more than a thousand decimal
places, which gives a perfect opportunity for tuning various
settings of the algorithm.

So it was without any expectations to find anything new that
the author performed the tests, which gave very satisfactory
results pertaining to the algorithm. However, it was
unsettling that some of the results were in disagreement
with the well known and long established facts such as
the Feigenbaum conjecture and the spectral properties of the doubling
or universality operator.

Thanks to the popularity of this field of mathematical and physical
sciences and to the Internet, most papers on the subject are readily
available. The present paper is a comparative study of the
spectral properties of the doubling operator
and a review of several works dealing with the original problem.
It is also an attempt to reconcile apparent contradictions and
to trace their origin.

Let us recapitulate briefly the setting of the problem.
It deals with the mapping of an interval onto itself $f\colon [a,b] \to [a,b]$,
where $f$ is a generic unimodal function. Here
unimodal means having only one extremum (maximum) on the interval $[a,b]$ (one-hump map),
and generic means that the function is smooth and the extremum is quadratic.
The function $f$ depends on one parameter.
The iterations of such maps can display an infinite cascade of
period doubling bifurcations as parameter changes.
The bifurcations occur when a stable solution $x_n$ to the equation
$x=f^{(n)}(x)$, $n=2^{k-1}$, $k \in \N$, loses stability, and two new stable
cycles are born, i.e., two solution $x_{2n}^{\pm}$ to the equation $x=f^{(2n)}(x)$.
It was shown in \cite{Fg1,Fg2}
with the help of renormalization involving rescaling and stretching of the iterated maps
that, as the period of cycles tends to infinity,
the sequence of bifurcations displays a universal character
independent of the initial function $f(x)$.
Asymptotically, the cascade of bifurcations possesses self-similarity with
the universal constants $\delta \approx 4.6$ in parameter space,
and $\alpha \approx -2.5$ in the phase space (on the interval).
These universal constants can be found from the period-doubling (or universality) equation:
\begin{eqnarray}
g(x)=T(g)(x) = g(g(g(1)x))/g(1), \quad x \in [-1,1].
\label{Un0}
\end{eqnarray}
Here the function $g(x)$ is the result of an infinite number of renormalizations
of iterations of the original mapping $f(x)$, and so it totally forgets its prehistory.
The constant $\alpha = 1/g(1)$; and the constant $\delta$ is determined from
the spectrum of the operator $dT(g)$, i.e., Fr\'echet derivative of the operator $T$
on the solution $g(x)$ to the universality equation (\ref{Un0}).
This solution cannot be found by iterations
of the operator $T$, since this operator is hyperbolic.
We recall that the Feigenbaum conjecture (in its modern form)
states that all the eigenvalues
of the operator $dT(g)$ except one lie within the unit circle,
the unstable eigenvalue being $\delta$. So the unstable manifold
at the fixed point $g(x)$ is one-dimensional.

A few remarks on the preceding paragraph. First, the equation (\ref{Un0})
is not a unique form of universality equation. In various papers, there are used
other forms of universality operator (see Sect.~2, 3). They possess the common
solution $g(x)$, which is an even analytical function in the neighborhood of
the interval $[-1,1]$. But these universality operators do not have the same spectrum,
and have different eigenfunctions for the same eigenvalues.
We will stick to the equation (\ref{Un0}) as canonical in this paper, and mark
the differences as they appear. Second, note the absence of the normalization
condition $g(0)=1$ usually imposed on the solution in the definition (\ref{Un0}).
The reason will be given in Sect.~4. Third, the renormalization used in
\cite{Fg1,Fg2} preserves the nature of the extremum of the original function $f$.
So the limit function must have the same type of extremum, and so the
equation (\ref{Un0}) must possess different solutions. This fact, of course, is
well known, and only mentioned here to avoid misunderstanding. We will deal
primarily with the function $g(x)$ having quadratic extremum, and discuss
other solutions in Sect.~4. The most important remark here is this:
the spectrum of an operator $dT(g)$ depends strongly on the functional
space where the operator is acting. For example, in the space $L_2[-1,1]$
of complex-valued functions integrable with square, the spectrum
of the operator $dT(g)$ is continuous and complex (Sect.~2).

Feigenbaum never mentioned specifically the functional space for
the operator $T$. In the framework of papers \cite{Fg1,Fg2} it is
hardly to be expected. But then, the function $g(x)$ had to be found,
and the most natural space for this is the space ${\cal E}$ of even bounded
analytical functions, since the function $g(x)$ must be smooth and even by construction.
Again, it was not made explicit, but the numerical algorithm described
in \cite[page 693]{Fg2} clearly uses discretization in the space ${\cal E}$
(or rather in its subset, see Sect.~4).
Since the finite dimensional approximation to the operator $dT(g)$
is obtained as a by-product of the Newton iterations scheme
used for numerical solution of (\ref{Un0}),
the spectrum found (numerically) in \cite{Fg2} corresponds to even eigenfunctions.
Hence, the Feigenbaum conjecture (Sect.~4).

In Sections 2, 3, we compare the spectrum of the operator $dT(g)$
and its various representations in different functional spaces.
We also discuss the most common mistakes made in various papers and monographs
in the analysis of the spectral problem for the universality
operator. Some mistakes are obvious as such, and some are the result
of misquoting or the wrong assumptions and peer pressure.

In Section 4, we solve numerically the spectral problem
for the universality operator in the Banach space ${\cal F}$ of bounded analytical
functions in an ellipse with the focal points $\pm 1$, with
the supremum norm, continuous on the closure of the ellipse.
Let us give a few reasons for this
choice. First, computer experiments revealed that the function $g(x)$
belongs to this space. A rigorous proof is still to be found,
despite some computer assisted efforts \cite{L1,L2}. But, as
``mathematics is an experimental science'' (V.I. Arnold),
we will consider this fact established.
The second reason is the fact that the functional space ${\cal F}$ admits
extremely good discretization with Chebyshev polynomials as a basis.
The coefficients of expansions of the functions in ${\cal F}$ in Chebyshev
series decrease exponentially \cite{Pa}.
This is why we choose an ellipse and not a disk
(see \cite[page 1264]{Ec}). Finally, and this is a purely physical argument
open for discussion: there is no reason to restrict
the space ${\cal F}$ to the space ${\cal E}$ of even analytical functions.
The function $g(x)$ forgets its prehistory and is even, but
all the pre-limit functions subject to renormalizations used in \cite{Fg1,Fg2}
still keep some information about the original function $f(x)$, which is unimodal,
and so the perturbations of the function $g(x)$ need not be even.
It is a part of universality that we need not impose some symmetry
on the function $f$ (as in the logistic map)
in order to obtain the function $g$.

\section{Explicit spectrum}

The notion of universality in dynamical systems
can now be found in almost every monograph remotely concerned with
chaotic dynamics.
An excellent exposition of the Feigenbaum universality can be found
in the book \cite[Chap.~7]{LL} aimed at physical scientists and engineers.
The book also illustrates how physical intuition fails
when simple mathematics is neglected. We will use this book as a typical
example.

The universality equation in \cite{LL} is given in the form
of rescaling equation (7.2.39) \cite[page 491]{LL}:
\begin{eqnarray}
g(x)=\alpha g(g(x/\alpha)) = T(g)(x), \quad g(0)=1.
\label{Un1}
\end{eqnarray}
Since the normalizing condition $g(0)=1$ is included in the definition,
it immediately follows that $\alpha = 1/g(1)$, and this equation
is identical with (\ref{Un0}). From the previous exposition in \cite{LL},
it also follows that the authors consider general maps, i.e., unimodal and generic
in the sense of Sect.~1, and so they implicitly operate in the space ${\cal F}$.

To investigate the stability of the fixed point $g(x)$, the authors
compute the linearized period-doubling operator
introducing a perturbation $g(x) + \varepsilon h(x)$ and,
linearizing, obtain the linear operator (G\^ateaux derivative) in the
form
\begin{eqnarray}
L(g)h(x) = \alpha \left( g'(g(x/\alpha)) h(x/\alpha) + h(g(x/\alpha))  \right).
\label{L}
\end{eqnarray}
Then the authors refer to Feigenbaum \cite{Fg1,Fg2} and claim (Feigenbaum conjecture)
that the spectrum
of the operator $L(g)$ has a single unstable (i.e., lying outside the unit circle)
eigenvalue $\delta \approx 4.669$ \cite[page 492]{LL}.
Unfortunately, both this claim and the linearized equation
(\ref{L}) itself are not true. So it is not clear how much of the following
physical argument in \cite{LL} will survive.

Let us compute the linearized period-doubling operator proceeding
exactly as described in \cite[page 491]{LL}, but keeping in mind that
$\alpha = 1/g(1)$, i.e., that $\alpha$ depends on the function $g$.

\noindent
{\bf Proposition 1.} {\it The formal G\^ateaux derivative of the operator $T$ defined in (\ref{Un0})
is given by the formula}
\begin{eqnarray}
dT(g)h(x) = L(g)h(x) + \alpha \left(g'(g(x/\alpha)) g'(x/\alpha) x - \alpha g(g(x/\alpha)) \right) h(1).
\label{dT1}
\end{eqnarray}

It seems that this easily verified formula (\ref{dT1}) for the operator $dT(g)$ was never computed.
An (almost) correct formula for the derivative was found in \cite{VS}, but for
another form of the universality operator (Sect.~3). The formula (\ref{dT1})
is applicable in any functional space where the G\^ateaux derivative of the operator $T$
coincides with the Fr\'echet derivative. It is certainly the case in the space ${\cal F}$.

\noindent
{\bf Proposition 2.} {\it The operator $T$ is compact (\cite[page 16]{VS}) in the space ${\cal F}$,
and the operator $dT(g)$ has the following spectrum
\begin{eqnarray}
S = [\lambda_1,\lambda_2,\dots]
 = [\alpha^2,\delta,\alpha,{1 \over \alpha},{1 \over \alpha^2},\lambda_6,{1 \over \alpha^3},\lambda_8,{1 \over \alpha^4},{1 \over \alpha^5},\dots],
\label{sp1}
\end{eqnarray}
where $|\lambda_i| > |\lambda_j|$, $i<j$.}

We will compute numerically the eigenvalues and the corresponding eigenfunctions of the
operator $dT(g)$ in Section 4. They coincide with $S$ in (\ref{sp1}). But now we
note that all the eigenvalues in $S$ where $\alpha$ is present are found explicitly
together with the corresponding eigenfunctions.

\noindent
{\bf Proposition 3.} {\it Let $k$ be any complex number except $1$. Then
$\lambda=\alpha^{1-k}$ is an eigenvalue of the formal
spectral problem $dT(g)h= \lambda h$ with the eigenfunction
\begin{eqnarray}
h(x)=g(x) - x g'(x) - g^k(x) + x^k g'(x).
\label{ak}
\end{eqnarray}
In addition, $\alpha^2$ is the eigenvalue with the eigenfunction
\begin{eqnarray}
h(x)=g(x) - x g'(x).
\label{a2}
\end{eqnarray}
}

Proof. If we differentiate the equation (\ref{Un1}) and put there $x=0$,
we obtain the identity $g'(1)=\alpha$. We use this, along with the
equation (\ref{Un1}) and its derivative, for simplifying substitutions.
We observe that if we put $x=0$ into the formal spectral problem
$dT(g)h(x)= \lambda h(x)$, then we derive the identity $(\alpha^2-\lambda) h(0) =0$.
Hence, for analytical functions $h$, either $\lambda=\alpha^2$, or $h(0)=0$.
The rest of the proof is a simple, although very bulky,
symbolic calculation better made on a computer \fullsquare

The spectral problem is formal until we specify the functional space
we are working with. In the space ${\cal F}$, obviously, $k=0,2,3,\dots$.
So we have found explicitly 7 out of the first 10 eigenvalues of the operator
$dT(g)$, and at least two of them lie outside of the unit circle.
This result is easily verified analytically (and numerically, Sect. 4)
and is in direct contradiction
with the Feigenbaum conjecture. So let us trace the origin of this
apparent paradox. But before we turn to the original paper \cite{Fg2},
where we hope to find an answer, we need the spectrum of the operator $L(g)$
for comparison.

\noindent
{\bf Proposition 4.} {\it The spectrum of the operator $L(g)$ in the space ${\cal F}$ is
\begin{eqnarray}
\tilde S = [\delta,\alpha,1,{1 \over \alpha},{1 \over \alpha^2},\lambda_6,{1 \over \alpha^3},\lambda_8,{1 \over \alpha^4},\dots],
\label{sp2}
\end{eqnarray}
where $\lambda_i$, $i=6,8,\dots$ are the same as in (\ref{sp1}).
The eigenvalues $\alpha^{1-k}$, $k=0,1,\dots$ in $\tilde S$
have the eigenfunctions
\begin{eqnarray}
h(x)=g^k(x)-x^k g'(x).
\label{ei2}
\end{eqnarray}
}

Proof. Numerically, it is demonstrated in Sect.~4. The part of the spectrum where
$\alpha$ is involved is found explicitly \fullsquare

We will return in Section 4 to all the spectral problems discussed in this
and the following Section.

\section{The problem with the spectrum}

Now we turn to the paper \cite{Fg2}, which is referenced in almost
every publication dealing with the Feigenbaum universality. To avoid confusion, we will
keep our notation $\alpha=1/g(1)<0$, which is common now
(Feigenbaum used $\alpha = -1/g(1) >0$, \cite[page 675]{Fg2}, \cite[page 73]{Fg3}),
and translate the corresponding formulas when needed.

Feigenbaum used a different form of the universality equation from what we use (\ref{Un0}).
It is given in the abstract of \cite{Fg2} as
\begin{eqnarray}
g(x)=\alpha g(g(-x/\alpha)) = T_2(g)(x).
\label{Un2}
\end{eqnarray}
The normalizing condition $g(0)=1$ is given later on in \cite{Fg2}.

Strangely, in the abstract of \cite{Fg2}, Feigenbaum gives
the linear operator ${\cal L}$, which coincides with (\ref{L})
on the function $g(x)$, since $g(x)$ is even.
The correct formula (assuming $\alpha={\rm const}$) should be
\begin{eqnarray}
L_2(g)h(x) = \alpha \left( g'(g(-x/\alpha)) h(-x/\alpha) + h(g(-x/\alpha))  \right),
\label{L2}
\end{eqnarray}
and the corresponding operator $dT_2(g)$ (correct Fr\'echet derivative) is
\begin{eqnarray}
\hspace*{-1.5cm}
dT_2(g)h(x) = L_2(g)h(x) - \alpha \left(g'(g(-x/\alpha)) g'(-x/\alpha) x + \alpha g(g(-x/\alpha)) \right) h(1).
\label{dT2}
\end{eqnarray}

Note that the formula (\ref{L2}) for the derivative of (\ref{Un2}) was found in
\cite[page 47, formula (42)]{Fg1}. The formulas (\ref{L2}) and (\ref{dT2})
can be simplified using the fact that $g$ is even and $g'$ is odd.
But this should be done after and not before the computation of the derivative
of the operator.
In addition, the function $h$ in these formulas need not be even,
so no simplifications there.

The spectral properties of the operators $L(g)$ and $L_2(g)$, and, respectively,
of the operators $dT(g)$ and $dT_2(g)$ are different in the space ${\cal F}$.
In the space ${\cal E}$, each pair of operators
possesses identical spectrum (Sect. 4).

Later on in \cite{Fg2}, Feigenbaum uses the operator $L(g)$ as the derivative of the
operator $T_2$ on $g(x)$, but periodically switches to $L_2(g)$
(see \cite[page 677, formula (28); page 682, 685]{Fg2}.

It is also not exactly clear, what Feigenbaum meant by his conjecture.
First, in the abstract of the paper \cite{Fg2}:
``${\cal L}$ possesses a unique eigenvalue in excess of 1.''
Then (we quote from \cite[page 687]{Fg2} using our notation and correcting a misprint):
``The spectrum of the operator $dT(g)$ is $\delta$ and $\alpha^{1-\rho}$,
$\rho \ge 1$, and, moreover, the spectrum is complete.''

We used here $dT(g)$ rather than $L(g)$, since here it was clearly meant
the derivative of the operator $T$.

The part about the spectrum being complete was refuted
numerically in many works, since
other eigenvalues were found (Sect.~4). In Proposition 2, they are
$\lambda_6$, $\lambda_8$, etc.

After numerical investigation of the spectral problem in \cite{Fg2},
Feigenbaum states his conjecture in the form \cite[page 694]{Fg2}:
``{\it $\delta$ is the solitary eigenvalue of $dT(g)$ greater than 1.}''

Note that all these conjectures imply that 1 is an
eigenvalue, and so they contradict the conjecture in its
modern interpretation. Although this difficulty is fixed
by the normalization $g(0)=1$, which simply means that we choose
one solution from the family of solutions, still, this eigenvalue is the product
of a wrong assumption. If the derivative $dT(g)$ was computed correctly,
the eigenvalue 1 would not appear (Sect.~4).

To complicate matters even more, Feigenbaum actually found
the eigenvalue $\alpha$, since $\rho=0$
perfectly fits the citation above, with the analytical eigenfunction
$1-g'(x)$ \cite[page 686]{Fg2}. This eigenvalue is not in excess of 1, since
$\alpha<0$, but $\alpha$ lies outside of the unit circle.

So let us draw a line here and try to explain these paradoxes.

First, Feigenbaum used the wrong linearization $L(g)$ instead of $L_2(g)$ of the universality operator.
In addition, both these linearizations are wrong, since
they assume $\alpha={\rm const}$ independent of the fixed point $g(x)$.
This assumption is later rejected in \cite[page 693, formula (80)]{Fg2},
when the variation of $\alpha$ is used together with the variation of $g(x)$.
The analysis of the spectrum is performed in some unspecified functional space, which is
clearly not a space of even functions, since some of the eigenfunctions (\ref{ei2})
are not even.

The second misunderstanding in \cite{Fg2} compounding the first
is the use of numerically obtained data in the same context
as analytically obtained eigenvalues and eigenfunctions.
These are two different sets of objects, since
the numerical algorithm described in \cite[page 693]{Fg2}
operates in a subset of the space ${\cal E}$
(see Sect.~4).
To unite the numerical and analytical data, we need the space ${\cal F}$
and correctly linearized operator $dT(g)$.

In the afterword to the paper \cite{Fg2}, Feigenbaum states that his spectral
conjecture was verified by Collet et al. We have no access to that paper (then in draft),
but in the subsequent publications of the same authors, the space of even functions
was postulated \cite[page 211, 212]{CEL}, \cite[page 427]{L1}, \cite[page 521]{L2}.

Now we consider how the spectral problem
for the universality operator was treated in several frequently cited papers
and in some books.

In the study \cite[page 4]{CEK}, the authors use the same notation as in this paper,
but consider the problem in a broader space of functions mapping
the interval $[-1,1]$ onto itself, i.e., the functions are not necessarily even.
The four assumptions, M1-M4, all agree with our conclusions so far, but then
the authors wrongly compute the derivative of the operator (\ref{Un0}) as
$L(g)$ (\ref{L}) and proceed with the analysis. In particular,
Lemma 1 in \cite{CEK} coincides with Proposition 4 here, so
the following assumptions M5, M6 \cite[page 5]{CEK} can be considered
as either true or wrong depending on what operator is taken for the derivative
of $T$. On the other hand, the authors found the eigenvalue 1, so
the solution $g$ to the equation $g=T(g)$ is either degenerate, or belongs
to a one parameter family of solutions (implicit function theorem).
Both facts are not true (Sect.~4).

In the paper \cite{Ec}, Eckmann gave some substantiation to the choice
of the space of even analytical functions, where $g(x)$
``is supposed to lie'', \cite[page 1264]{Ec}. His space is similar
to the space ${\cal E}$, except it is defined on a disk, not an ellipse.
However, the properties P1-P3 (including the Feigenbaum conjecture)
hold there only with an additional stipulation (see Sect.~4).

Feigenbaum renormalizations preserve the property of the function $f$
being symmetric with respect to its hump, so the choice
of the functional space of even functions is justified for such maps
(logistic map, for example). But the Feigenbaum universality is now understood in
a broader sense (see \cite[Chap.~7]{LL}), meaning the functions $f$ need not be symmetric.
This confusion of notions leads to many erroneous statements on
the dimension of the unstable manifold at $g(x)$. For example,
in the paper \cite[page 425]{Su}, the author refers to
Lanford's computer-assisted proof, but explains his results
in a general space of analytic functions;
in the book \cite{Sh}, analytical unimodal maps are considered,
so Proposition 2 in \cite[page 191]{Sh} and its corollaries are not true;
in the book \cite{Ar}, the Feigenbaum universality is explained
on a typical example $f(x)=A x \exp(-x)$ \cite[page 338, 339]{Ar},
but the doubling operator $J$, identical to $T$ (\ref{Un0}),
is defined on even functions with some restrictions \cite[page 340]{Ar},
and the Feigenbaum conjecture is formulated in an unspecified
functional space.

Although we are not concerned with proofs of the Feigenbaum conjecture in this paper,
some of the works on the subject deserve a special attention,
since they apparently disagree with our results stated above.

Lanford is reputed to have given the first of
the computer-assisted proofs of the Feigenbaum conjecture.
His proof seemingly contradicts our conclusions,
but only if the results are taken out of the context.
In the paper \cite{L1}, he introduces the space ${\cal M}$ of continuously differentiable
even mappings $\psi$ of the interval $[-1,1]$ into itself such that
$\psi(0)=1$ (among other things) \cite[page 427]{L1}. But the condition $\psi(0)=1$
makes ${\cal M}$ a set, not a space, since
functions cannot be added or multiplied by a constant in ${\cal M}$.
Further \cite[page 428]{L1}, he introduces a Banach space ${\cal B}$ of bounded
even analytic functions on a set $\{z \in \CC \colon |z^2-1| <2.5\}$
equipped with the supremum norm,
and its subspace ${\cal B}_0$ of functions vanishing to second order at 0.
Theorem 3 on hyperbolicity of $dT(g)$ \cite[page 428]{L1} is
formulated in the space ${\cal B}_0$, where it is not true,
since the functions in this space do not satisfy the universality equation.
It is, probably, a misprint, since Theorem 3 is true in the set (or an affine
space) ${\cal B}_1={\cal B}_0+1$ (see Sect.~4).
Further \cite[page 429]{L1}, Lanford introduces the expansion
$\psi(x) =  1 -  x^2 h(x^2)$ corresponding to the set ${\cal B}_1$,
which was used in many papers implicitly.

In the paper \cite{EcW}, where another computer assisted proof
of the Feigenbaum conjectures is given, the word ``even'' is not mentioned
even once. However, even functions are implied \cite[Theorem 2.2]{EcW}.
It is also the case in \cite[page 396]{Ep} and many other papers.

To the best of the author's knowledge, there is a unique paper \cite{VS}
where the correct formula for the derivative of the doubling operator
was found (but for the wrong operator).
The authors consider generic unimodal maps as defined in Sect.~1 \cite[page 14]{VS}
(we quote the Russian edition), and the Feigenbaum conjecture is formulated
in its modern form without reference to even maps.
The doubling operator is defined \cite[page 13]{VS} as
\begin{eqnarray}
T_3(g)(x) = -a g(g(x/a)), \; a = - {g(0) \over g(g(0))}, \; g(0) ={\rm const}, \quad x \in [-1,1].
\label{Un3}
\end{eqnarray}
Here we substitute $a$ for $\alpha$ to avoid confusion. If $g(0)=1$, then $a=-\alpha$.
The authors compute the correct derivative \cite[page 16]{VS}, but for the operator
\begin{eqnarray}
T_4(g)(x) = -a g(g(-x/a)),
\label{Un4}
\end{eqnarray}
which is not the same as (\ref{Un3}) for analytic functions.
The analysis of spectral properties of the operator $dT_4(g)$ in \cite{VS}
is very similar to that in the present paper, although it is more difficult
due to a more complicated form of the
doubling operator. The authors have found the eigenvalue 1, and $-a=\alpha \approx -2.5$,
as well as other powers of $\alpha$, except $\alpha^2$.
These results contradict the Feigenbaum conjecture stated earlier in \cite{VS}. So
the authors have tried to dismiss unwanted eigenvalues on the following grounds \cite[page 17]{VS}.
First, they are not relevant to the universality, since they are linked to
coordinate transformations (i.e., not to the parameter space).
Not many people would subscribe to this point of view today, since $\alpha$
is considered now on a par with $\delta$ as a universal constant.
For example, in \cite[page 488]{LL}, it is explained that
both $\alpha$ and $\alpha^2$ play a part in the rescaling of periodic
solutions. The second argument the authors use to conform to the Feigenbaum conjecture
is (a) -- the eigenvalue 1 is eliminated by the condition $g(0) ={\rm const}$;
and (b) -- the eigenvalue $\alpha$ is eliminated by the condition $g'(0) =0$.
The condition (a) means that we choose one solution from a family,
so the eigenvalue 1 is simply ignored;
and the condition (b) was not imposed in the statement of the problem, and anyway,
it follows from the universality equation, i.e., $g'(0) (\alpha-1) =0$ follows from (\ref{Un0}), and
similarly for (\ref{Un4}).
The property $g'(0) =0$ of the solution $g(x)$ to (\ref{Un0}) or (\ref{Un4}) is a result
of an infinite number of renormalizations. But perturbations of the solution
need not conform to this restriction.
In addition, this projection does not explain what to do with other powers of $\alpha$ present in
the spectrum in both spaces ${\cal F}$ and ${\cal E}$ (Sect.~4).
Further, the authors give incorrect form of the doubling operator \cite[page 19, formula 4.1]{VS}
with $\alpha=g(1)$, but this is clearly a misprint.

In some papers, the derivative of the operator (\ref{Un0}) is computed
incorrectly, but then never used; so the mistake is not revealed
\cite[page L713S]{CCR}. And the use of the space of even functions
can only be deduced by a dedicated reader. In \cite[page L713S]{CCR},
it was only indicated as ``Lanford's expansion'' of the function $g(x)$.

We conclude this survey with two works devoted to precise
computation of the Feigenbaum constants.

In the paper \cite{Br1}, Briggs
uses the same notation and the same operator (\ref{Un0}) as we used in this
paper \cite[formula (5)]{Br1}. Numerical algorithm is
similar to that used by Feigenbaum and the most authors \cite[page 437]{Br1}, i.e.,
it operates in the subset of the space of even functions.
The Feigenbaum conjecture is formulated for historical reference;
then the wrong ``local linearization of $T$'' \cite[formula (8)]{Br1}
is obtained by ``simple calculation''. In fact, this local linearization coincides
with that of Feigenbaum in \cite[page 47, formula (42)]{Fg1}, where it is
found for the positive $\alpha=-1/g(1)$.
Fortunately, this ``local linearization'' was never used in \cite{Br1}.

In his PhD thesis \cite{Br2}, Briggs uses the same notation as in \cite{Br1}
(see \cite[formula (1.5)]{Br2}). But then, the derivative $DT_g$ of the
operator $T$ is upgraded to include the dependence of $\alpha$ on the solution $g(x)$
\cite[page 5]{Br2}. This new formula for the derivative $DT_g$ is remarkably similar
to that found in \cite[page 16]{VS} for the different operator (\ref{Un4}).
Then, \cite[page 12]{Br2}, the ``local linearization of $T$''
is found again by ``simple calculation'' as in \cite{Br1},
but this time with the correct sign of $\alpha$. Apparently, Briggs is familiar with
eigenvalues which do not comply with the Feigenbaum conjecture, but he explains them
as ``extra eigenvalues'' introduced by a finite-dimensional approximation \cite[page 22]{Br2}.
Briggs recommends to select the good eigenvalues, which are readily identified,
and discard the bad ones.

In the next Section, we will not follow this advice.

\section{Numerical analysis of the spectral problems}

In this Section, we compute the spectrum for all spectral problems
mentioned in previous Sections in different functional spaces.
We will also use different algorithms including that described by
Feigenbaum in \cite[page 693]{Fg2}, which was used (with various modifications)
by Lanford \cite{L1}, Briggs \cite{Br1}, and many other researches.

First, we describe an algorithm based on the use of Chebyshev polynomials
as a basis in the space ${\cal F}$.

The solution $g(x)$ to the equation (\ref{Un0}) is approximated
by the polynomial
\begin{eqnarray}
g(x) = \sum_{i=1}^n g(x_i) p_{ni}(x), \quad x_i=\cos {(2i-1)\pi \over 2n}, \quad i=1,\dots,n,
\label{Pg}
\end{eqnarray}
where $x_i$ are Chebyshev roots, and
\begin{eqnarray}
p_{ni}(x) = {T_n(x) \over (x-x_i) T'_n(x_i)}, \quad i=1,\dots,n
%\label{pnx}
\end{eqnarray}
are Chebyshev fundamental polynomials of Lagrange interpolation.
We will use the notation $g(x)$ both for the analytic solution to (\ref{Un0})
and for its polynomial approximation (\ref{Pg}) (and others),
but this will not lead to confusion.

The equation (\ref{Un0}) is rewritten as $\Phi(g)=g-T(g)=0$, and the solution
is found by the Newton iterations
$$
g_{k+1} = g_k - A_k^{-1} \Phi(g_k), \quad k=0,1,\dots,
$$
where $g_k$ is the $k$-th approximation to the solution $g$;
$A_k=d \Phi(g_k)$ is the Jacobian matrix at $g_k$.
The iterations are done until the polynomial $g_{k+1} - g_k$ (evaluated at
the nodes $x_i$) is zero in the sup norm
within the round-off error.

After the final iteration, we have found
the polynomial solution $g$ represented
by the values $\{g(x_i), i=1,\dots,n\}$, and
the matrix $A=A(g)$ on the solution.
The matrix $I-A$ is an approximation to the derivative $dT(g)$ (\ref{dT1}),
where $I$ is the unit matrix. Then we compute the spectrum
of the matrix $I-A$ by standard linear algebra subroutines.

In the course of these computations, we need to evaluate the polynomial
$g(x)$ at the points that are not Chebyshev roots.
This can be done very efficiently, if the function $g(x)$ is expanded
in Chebyshev series
\begin{eqnarray}
g(x) = a_0/2 + \sum_{k=1}^{n-1} a_k T_k(x), \quad T_k(x)=\cos(k \arccos x),
\label{Chg}
\end{eqnarray}
where $T_k(x)$ are Chebyshev polynomials.
The coefficients $a_k$ are found by the discrete Fourier-Chebyshev
transform. This operation is stable and does not accumulate the round-off
errors \cite{Pa}. The evaluation of $g$ is done with the series (\ref{Chg})
using the recurrence relations for Chebyshev polynomials.
These operations are also stable \cite{Pa}.

The Fourier-Chebyshev transform also provides a built in precision control,
since the coefficients $\{a_k, k=1,\dots,n\}$ must decrease exponentially.
This can be seen on a plot of $\log(1/|a_k|)$ versus $k$, $k=1,\dots,n$.

The elements of Jacobian matrix $A$ are approximated by finite differences.
It needs not be done with high precision for Newton iterations to
converge quadratically. Only after the final iteration this
matrix needs to be evaluated with maximal precision,
since it is used for the approximation of the spectral problem.

We have also used an alternative way to approximate the operator $dT(g)$ (\ref{dT1}).
If the doubling operator $T$ (\ref{Un0}) is applied to a polynomial $p$ of
an order $m$, then $T(p)$ is a polynomial of the order $m^2$. Since $T(p)$
can be restored by its values at $m^2+1$ points, the same is true for
the derivative $dT(p)$. So if we take the dimension $n$ of the projection
such that $n \ge m^2+1$, then we can compute the operator $dT(p)$
exactly, i.e., in the same sense as Gauss quadratures are
exact on polynomials up to a certain order.

The finite difference approximations are proved to be faster, but
slightly less accurate.

This algorithm takes about as many lines in a computer language
as it took to describe it. For general analytic functions in ${\cal F}$,
it is also one of the most efficient. It follows from the approximative
properties of the Chebyshev series (\ref{Chg}) and
asymptotically optimal distribution of Chebyshev nodes (see \cite{Pa}).
However, for the same accuracy of the solution $g$,
this algorithm takes about 4 times more memory and 8 times
more CPU time than the algorithms that use the symmetry
of the solution $g(x)$. This is probably why it was never used before.
Recently, Chebyshev series representation of $g(x)$ was
used in \cite{Ma}, but on the interval $[0,1]$, i.e., for even functions.

We are not about to break any records in the number of
digits of the universal constants. The original plan was to test
the algorithm, so we fix the number of nodes on the interval $[-1,1]$ as $n=32$,
and we fix the floating point arithmetic at 64 digits. The software
for such computations is available as an open source (see \cite{Bai}).
The chosen precision is equivalent to working with infinite number
of digits, since the round-off errors can be neglected in comparison
with the errors of the approximation. We should mention that
all computations were verified with different settings
(more/less digits/nodes, and different linear algebra routines
for solution of the spectral problems),
but they gave similar results and not reported here.

The constant $\alpha=1/g(1)$ is found with the accuracy $0.5 \times 10^{-22}$,
which is a small number in comparison with
the last coefficient ($\approx 0.2454065396 \times 10^{-13}$)
at $x^{30}$ in the Taylor expansion of $g(x)$. The reason for this
is the value of the last Chebyshev coefficient at $T_{30}(x)$
in the expansion (\ref{Chg}), which is $0.4571053006 \times 10^{-22}$.
The constant $\delta$ is found with 22 correct decimal places.
We stress that no normalization needs to be imposed on the solution
$g(x)$. The Newton iterations converge quadratically,
provided a good initial approximation is taken, and the
solution is found uniquely in the space ${\cal F}$.

In Table 1, we cite the first 11 eigenvalues of
the operator $dT(g)$ computed as described above. They correspond
to the spectrum $S$ in Proposition 2 to
the indicated accuracy, which was estimated by comparison
with the values of $\alpha$ and $\delta$ found in \cite{Br1}.
We cite here only 10 decimal places and can send the computed values on demand.

\centerline{{\bf Table 1.} First 11 eigenvalues of the spectrum $S$.}
$$
\begin{tabular}{lr|c|c}
%\hline
$\lambda_1=$      &\fm $ 6.264547831$ & $\alpha^2$      & $0.7 \times 10^{-21}$   \\
$\lambda_2=$      &\fm $ 4.669201609$ & $\delta$        & $0.2 \times 10^{-21}$   \\
$\lambda_3=$      &    $-2.502907875$ & $\alpha$        & $0.2 \times 10^{-21}$   \\
$\lambda_4=$      &    $-0.399535280$ & $\alpha^{-1}$   & $0.5 \times 10^{-21}$   \\
$\lambda_5=$      &\fm $ 0.159628440$ & $\alpha^{-2}$   & $0.4 \times 10^{-18}$   \\
$\lambda_6=$      &    $-0.123652712$ &                 &                         \\
$\lambda_7=$      &    $-0.063777193$ & $\alpha^{-3}$   & $0.3 \times 10^{-18}$   \\
$\lambda_8=$      &    $-0.057307021$ &                 &                         \\
$\lambda_9=$      &\fm $ 0.025481238$ & $\alpha^{-4}$   & $0.1 \times 10^{-12}$   \\
$\lambda_{10}=$   &    $-0.010180653$ & $\alpha^{-5}$   & $0.9 \times 10^{-17}$   \\
$\lambda_{11}=$   &    $-0.010145805$ &                 &                         \\
\end{tabular}
$$

The eigenvalues that correspond to the powers of $\alpha$ in Table 1 have
the eigenfunctions given in Proposition 3 (after the normalization).
The eigenvalues $\lambda_2=\delta$, $\lambda_6$, $\lambda_8$, $\lambda_{11}$, etc.,
that are not related to the powers of $\alpha$ (at least not in an obvious way),
have even eigenfunctions.

The spectrum of the operator $L(g)$ (\ref{L}), that
was frequently mistaken for the derivative of the operator $T$,
is given in Proposition 4. To verify it numerically,
we only need to fix the value of $\alpha$ in the program,
after $g(x)$ and $\alpha$ are already found.
The eigenvalue 1 indicates that there is a one-parameter family
of solutions. We will find the family below for another problem.

The spectrum of the operator $dT_2(g)$ (\ref{dT2}) in the space ${\cal F}$
is
$$
S_2 = [\alpha^2,\delta,-\alpha,-{1 \over \alpha},{1 \over \alpha^2},\lambda_6,-{1 \over \alpha^3},\lambda_8,{1 \over \alpha^4},-{1 \over \alpha^5},\dots],
$$
i.e., $\alpha$ in $S$ (\ref{sp1}) is replaced with $-\alpha$. The same is true
for the operator $L_2(g)$ (\ref{L2}), i.e., it is the spectrum $\tilde S$
of the operator $L(g)$ (\ref{L}) with the substitution $\alpha \to -\alpha$.
However, only even eigenfunctions in Propositions 3, 4
are preserved, i.e., (\ref{a2}) for $\alpha^2$, and
(\ref{ak}), (\ref{ei2}) for odd $k$. For even $k$, we could not find explicit formulas.

Now we compute the spectrum of the operators used in \cite{VS}, i.e., $T_4(g)$ (\ref{Un4}), and
$T_3(g)$ (\ref{Un3}). We recall that in \cite{VS}, the doubling operator was defined
as $T_3(g)$, but the derivative was computed for the operator $T_4(g)$.
Since both operators have 1 as an eigenvalue, the Newton iterations do not converge
(i.e., there is a family of solutions). So we used $g(x)$ found earlier, which
satisfies all universality equations.
The spectrum of the operator $dT_4(g)$ is $\tilde S$ (\ref{sp2}), i.e.,
it coincides with the spectrum of the operator $L(g)$.
However, the explicit eigenfunctions corresponding to $\lambda=\alpha^{1-k}$, $k=0,2,3,\dots$
are the same as for the operator $dT(g)$ in Proposition 3;
and the eigenvalue 1 has the eigenfunction $g(x) - x g'(x)$ (as was found in \cite{VS}).
The eigenvalue $\alpha^2$ is missing in these problems.

It turns out that the spectra of $T_4(g)$ and $T_3(g)$ stand in the same relationship
as the spectra of $T(g)$ and $T_2(g)$, i.e., the spectrum of $dT_3(g)$ is obtained
from the spectrum of $dT_4(g)$ by the substitution $\alpha \to -\alpha$,
and half of the explicit eigenfunctions (for even $k$) could not be recovered.

The one-parameter family of solutions to the equations (\ref{Un1}) and (\ref{Un4}) can be found
if we take $g(0)$ as a parameter on the family and fix it in the procedure.
Numerical solutions that we found correspond to the family $g_\mu(x)=\mu g(x/\mu)$, $\mu \in \R$.
The value $\alpha$ and the spectrum are preserved on the family;
however, only one explicit eigenfunction $g_\mu(x) - x g_\mu'(x)$ is left
for the eigenvalue 1 in each problem.

Thus, the equation $y(x) = \beta y(y(x/\beta))$ has a family
of solutions only for a discrete set of values $\beta$. One of them is $\beta=\alpha \approx -2.5$,
another is $\beta=1$ for $y(x) \equiv 1$.
For each family, there is a solution $y_0(x)$ on the family for which $\beta=1/y_0(1)$.
For $\beta=\alpha$, $y_0(x)=g(x)$ (the solution to the equation (\ref{Un0}));
for $\beta=1$, $y_0(x) \equiv 1$, and the family itself is $y(x) = {\rm const}$ with the
spectrum $1,0,0,\dots$
Since the spectrum is preserved on each family of solutions, the families cannot intersect.

Other families can be found for different types of extrema of the unimodal
solution to the equation (\ref{Un0}), i.e., $g(x)-1=O(x^{2k})$, $k=2,3,\dots$
We found the solution for $k=2$ (with $n=70$)
\begin{equation}
\begin{array}{l}
g(x)=1 - 1.834107907\, x^4 + 0.012962226 \, x^8 + 0.311901736 \, x^{12}  \cr
\phantom{g(x)=} -0.062014622 \, x^{16}  -0.037539249 \, x^{20} + 0.017665496 \, x^{24} + \dots
\end{array}
\label{g_2} 
\end{equation}
for which $1/g(1) = \alpha_2 \approx -1.690302971$
with the accuracy $0.1 \times 10^{-19}$, as compared to \cite{Br1}.
The constant $\delta_2 \approx 7.284686217$ is found with 19 correct decimal places.
The spectrum of the operator $dT$ on the solution (\ref{g_2}) is
$$
[\gamma_1,\gamma_2,\dots] =
[\alpha_2^4,\delta_2,\alpha_2^3,\alpha_2^2,\alpha_2,{1 \over \alpha_2},{1 \over \alpha_2^2},\gamma_8,\gamma_9,{1 \over \alpha_2^3},\dots],
$$
where $|\gamma_i| > |\gamma_j|$, $i<j$, and $\gamma_8 \approx 0.291838408$, 
$\gamma_9 \approx -0.255664558$. Proposition 3 holds here formally except for
the eigenfunction (\ref{a2}), which corresponds now to the eigenvalue $\alpha_2^4$.

Now we turn to different functional spaces, and, to make it more demonstrative,
we will use a different quasi-numerical algorithm.

We will consider various Taylor expansions of the solution $g(x)$ to the equation (\ref{Un0}).
The coefficients of these series are found exactly (symbolically) in rational arithmetic
by the symbolic values of the polynomial $g(x)$ at the chosen rational nodes.
Symbolic approach avoids the floating point arithmetic at this very crucial stage,
since the corresponding linear systems of equations are very ill-conditioned.
In this way, we obtain an analog of the Fourier-Chebyshev transform
for arbitrary distributed rational nodes.
Thus, the floating point arithmetic is only used for the evaluation
of polynomials.

Let us verify the Feigenbaum conjecture for the equation (\ref{Un0}) using this algorithm and the
Lanford's expansion $g(x) = 1 - x^2 y(x^2)$. As it was mentioned,
this substitution is frequently used for the numerical
solution of the equation (\ref{Un0}) (including the paper \cite[page 693]{Fg2}).

The polynomial $g(x)$ is expanded in the Taylor series
\begin{eqnarray}
g(x)=1+a_1 x^2 + a_2 x^4 + \dots + a_m x^{2m},
\label{Br}
\end{eqnarray}
where $m$ is the dimension of the approximation and the number of nodes
taken on the interval $[0,1]$. This set of functions is not a space,
but a subset in the space ${\cal E}$ of even analytical functions.
In addition, true eigenfunctions do not belong to this set, since
they all (except one) vanish at the origin (see Proposition 3).

We take $m=15$, and choose the nodes $x_i=i/m$, $i=1,2,\dots,m$.
Then we solve symbolically the linear system $\{g(x_i)=g_i\}$, $i=1,2,\dots,m$
with respect to the coefficients $a_k$, $k=1,2,\dots,m$ of the
Taylor expansion (\ref{Br}). Then we evaluate this exact
solution as needed in floating point arithmetic on different
sets of values $g(x_i)$, $i=1,2,\dots,m$.
The Newton iterations are done as described above for Chebyshev nodes,
and the spectral problem is solved after the final iteration
for the obtained matrix $I-A$. Thus we find the spectrum
$$
S_3 = [\delta,{1 \over \alpha^2},\lambda_6,\lambda_8,{1 \over \alpha^4},\dots],
$$
i.e., we recovered the Feigenbaum conjecture. Only those eigenvalues of the spectrum $S$ (\ref{sp1})
are left in $S_3$ that correspond to even eigenfunctions (except for $\alpha^2$).
The constant $\delta$ is found with 19 correct decimal places,
although the power of the polynomial solutions is the same, i.e., 30.
This is due to the poor choice of the nodes, compensated only
by the rational arithmetic.

In many papers, only the space of even analytical functions is defined,
and the Lanford's expansion (\ref{Br}) is not stipulated
(see, for example, \cite[page 1264]{Ec}). This makes the Feigenbaum conjecture
not true. To demonstrate this, we take the expansion
\begin{eqnarray}
g(x)=a_0+a_1 x^2 + a_2 x^4 + \dots + a_m x^{2m},
\label{Ec}
\end{eqnarray}
and proceed as described above, but for $x_i=i/m$, $i=0,1,2,\dots,m$.
We obtain the spectrum $S_3$ plus the missing eigenvalue $\alpha^2 \approx 6.26$.

In both cases (\ref{Br}) and (\ref{Ec}), the same solution $g(x)$ is obtained,
and the Newton iterations converge quadratically (as the spectrum indicates they should).

Now we demonstrate that even functions are not necessary for the Feigenbaum
conjecture to be fulfilled.

First, we take the expansion
\begin{eqnarray}
g(x)=a_0+a_1 x + a_2 x^2 + \dots + a_m x^m
\label{V1}
\end{eqnarray}
on the interval $[-1,1]$. We take the Chebyshev nodes and approximate them
as rational numbers (with small denominators).
Then we solve symbolically the linear system $\{g(x_i)=g_i\}$, $i=0,1,2,\dots,m$
and proceed as described above for the quasi-numerical algorithm.
This is another projection of the space ${\cal F}$ on a
finite dimensional one. As expected, we duplicated the results
obtained with the Chebyshev approximation (\ref{Pg}) and
obtained the spectrum $S$ (\ref{sp1}).

Now we fix $a_0=1$ in the expansion (\ref{V1}), decrease the dimension by 1,
and repeat the process. This will kill the eigenvalue $\alpha^2$.
If we keep $a_0$ arbitrary and fix the coefficient $a_1=0$,
then we kill the eigenvalue $\alpha$, but $\alpha^2$ is still present.
Finally, if we fix both $a_0=1$ and $a_1=0$ and repeat the process,
then we kill both eigenvalues $\alpha^2$ and  $\alpha$ and
recover the Feigenbaum conjecture.

\section{Conclusion}

In the paper \cite{CEL}, the condition P3), i.e., the functions
being even, is imposed ``mostly for convenience; it simplifies matters
and is satisfied by the $\psi$'s we are able to analyze in detail'' \cite[page 211, 212]{CEL}.
A rhetorical question is: how much our convenience and ability
to analyze something in detail are related to the physical
relevance of such an analysis? We do not pretend to know the answer
to this question with respect to the Feigenbaum universality.
However, in other problems, for example, bifurcations of periodic solutions
in a dynamical system, there is no reason to restrict the analysis
to symmetric functions if the solutions in question possess the symmetry.
On the contrary, the loss of the symmetry is one of the possible bifurcations
(see \cite{VV}).

These considerations lead us to believe that the Feigenbaum conjecture
in its present form is a numerical artifact. It is not clear,
why so much effort was spent on elimination of both additional
unstable eigenvalues $\alpha$ and $\alpha^2$, since they
both play a part in the rescaling of periodic
solutions (see \cite[page 488]{LL}).

\end{document}